\documentclass{amsart}
\usepackage{graphicx} 

\usepackage{amsmath,amssymb}
\usepackage{indentfirst}
\usepackage{hyperref}
\usepackage{mathrsfs}
\usepackage{amsthm}
\usepackage{mathtools}
\usepackage{enumitem}
\usepackage{titlesec}
\usepackage{comment}
\allowdisplaybreaks[1]

\newtheoremstyle{spacedthm}
  {\topsep}
  {18pt plus 4pt minus 2pt}
  {\itshape}
  {}
  {\bfseries}
  {.}
  {.5em}
  {}

\newtheorem{theorem}{Theorem}
\newtheorem{lemma}[theorem]{Lemma}

\newtheorem{proposition}[theorem]{Proposition}

\makeatletter
\renewenvironment{proof}[1][\proofname]
{\par
  \pushQED{\qed}%
  \normalfont
  \topsep6\p@\@plus6\p@\relax
  \trivlist
  \item[\hskip\labelsep\itshape #1\@addpunct{.}]
}
{\popQED\endtrivlist\@endpefalse}
\makeatother
\title{Holomorphic planes in $M\times \mathbb{T}^2$ from weakly exact Lagrangians in $M$}
\author{Urs Fuchs and Yuchen Yang}
\address{Department of Mathematics\\ University of Aachen\\ 52062 Aachen}

\email{fuchs@mathga.rwth-aachen.de, yuchen.yang@rwth-aachen.de}
\date{November 2026}

\titlespacing{\subsection}
  {0pt}
  {5ex}
  {2ex}

\begin{document}
\begin{abstract}
    Let \((M,\omega_M)\) be closed symplectic manifold containing a closed Lagrangian  submanifold \(L\) with $\omega_M$ vanishing on
    \(\pi_2(M,L)\). We prove that there is a nonconstant \(J\)-holomorphic plane $u:\mathbb{C}\rightarrow  M\times \mathbb{T}^2$ for any almost complex structure $ J$ on \(M\times \mathbb{T}^2\) tamed by \(\omega_M\oplus\omega_{\text{std}}\).
\end{abstract}
\maketitle

\addcontentsline{toc}{section}{Introduction}

In a recent article \cite{cattalani2025ahlfors_classic}, Cattalani constructed nonconstant $J$-holomorphic planes in a closed symplectic manifold $M\times \mathbb{T}^2$ with a split symplectic form taming $J$. His construction relies on the following assumptions: (i) $M$ is symplectically aspherical, (ii) $M$ retracts to a closed Lagrangian $L\subset M$ and (iii) the homological Dehn function of $M$ grows at most quadratically. Assumptions (i) and (ii) imply that $L$ is weakly exact, i.e. that $\omega$ vanishes on $\pi_2(M,L)$.

The goal of the present article is to show that such an existence result still holds if $L$ is just assumed to be weakly exact:

\begin{theorem}\label{main; theorem} Let \((M,\omega_M)\) be closed symplectic manifold and $L\subset M$ a closed Lagrangian submanifold such that $\omega_M$ vanishes on
    \(\pi_2(M,L)\). \\Then for any almost complex structure $ J$ on \(M\times \mathbb{T}^2\) tamed by \(\omega_M\oplus\omega_{\text{std}}\) there is a nonconstant \(J\)-holomorphic plane $u:\mathbb{C}\rightarrow  M\times \mathbb{T}^2$.
\end{theorem}

 To prove this Theorem we generally follow Cattalani's proof: We first use the Lagrangian $L\subset M$ to obtain by Theorem \ref{Lem; existence of holomorphic disk} a family of $J$-holomorphic disks $u_R:(\mathbb{D},\partial\mathbb{D})\rightarrow  (M\times \mathbb{R}^2, L\times \partial\mathbb{D}_R)$. Then we extract subdisks $\Delta_R\subset\mathbb{D}$  with $\text{length}(u_R(\partial \Delta_R))=o(\text{area}(u_R(\Delta_R)))$ such that Duval's Theorem \ref{thm; Duval} applies to the restrictions $u_R\vert_{\Delta_R}$ to obtain nonconstant $J$-holomorphic planes in $M\times \mathbb{T}^2$.\\  Cattalani chose the boundary $\partial \Delta_R$ such that $u_R(\partial\Delta_R)$ lies within a 1-neighborhood of $L\times \partial\mathbb{D}_R$ such that $\text{length}(u_R(\partial\Delta_R))=o(\text{area}(u_R(\Delta_R)))$. Instead, we use the area bound on $u_R(\mathbb{D})$ from Proposition \ref{remark; area of D with the pull back metric} to choose $\partial \Delta_R$ such that $u_R(\partial \Delta_R)$ lies within distance at most $R/2$ of $\partial \mathbb{D}_R$ (when projected to $\mathbb{R}^2$) and such that $\text{length}(u_R(\partial \Delta_R))= O(R)$. Then $\text{area}(u_R(\Delta_R))$ still grows quadratically in $R$ and we can avoid assumptions beyond the existence of a weakly exact Lagrangian $L$.\\

Above and throughout the paper, statements about area and length concern the parametrized area and length counted with multiplicity, obtained by pulling back a suitable Riemannian metric on the target space. In the following, if $u:\Sigma\rightarrow X$ is a smooth map from a surface $\Sigma$ into a Riemannian manifold $(X,h)$, we will denote them by
$\text{area}_{u^*h}(\Sigma)$ resp. $\text{length}_{u^*h}(\partial \Sigma)$.\\ We denote $\mathbb{D}_R:=\{x\in \mathbb{R}^2\, \vert |x|\leq R\}\subset \mathbb{R}^2$ be the round disk of radius $R$ in $\mathbb{R}^2$ and we also write simply $\mathbb{D}$ for $\mathbb{D}_1$.\\
Let $g$ be the Riemannian metric on $M\times \mathbb{T}^2$ defined by $$g(v,v'):=\frac{1}{2}(\omega(v,Jv')+\omega(v',Jv)).$$
Let $\pi:M\times \mathbb{R}^2 \rightarrow M\times \mathbb{T}^2$ be the covering map induced by the universal covering map $\mathbb{R}^2\rightarrow \mathbb{T}^2$.
We denote by $\tilde{J},\tilde{\omega},\tilde{g}$ the pullback of $J, \omega$ resp. $g$ via $\pi$ to $M\times \mathbb{R}^2$.\\

The following result is a a special case of a result by Gromov \cite[Theorem 0.4.A]{Gromov1985}, here for $V=M\times \mathbb{R}^2$ and the Lagrangian submanifold $W=L\times \partial \mathbb{D}_R$. Such a result is also essential in Cattalani's proof, see \cite[Lemma~2.25]{cattalani2025ahlfors_classic} for the corresponding statement and an outline of a proof.
\begin{theorem}
\label{Lem; existence of holomorphic disk}
     In the situation above, there is for every $R>0$ a \(\widetilde{J}\)-holomorphic disk
    \(u_R:(\mathbb{D},\partial\mathbb{D}) \to (M\times \mathbb{R}^2,L\times \partial\mathbb{D}_R)\)
    such that
    \(\operatorname{deg}(p_2\circ u_R|_{\partial\mathbb{D}}:\partial \mathbb{D}\rightarrow \partial \mathbb{D}_R)=1\).\\
    Here $p_2:M\times \mathbb{R}^2\to \mathbb{R}^2$ denotes the projection to the $\mathbb{R}^2$ factor. \hfill \qed
\end{theorem}

\noindent
\begin{proposition}
    \label{remark; area of D with the pull back metric}
Under the assumptions of Theorem~\ref{Lem; existence of holomorphic disk} we have $\text{area}_{u_R^*\tilde{g}}(\mathbb{D})=\pi R^2$.
\end{proposition}

\begin{proof}[Proof of Proposition~\ref{remark; area of D with the pull back metric}]
The symplectic form $\tilde{\omega}:=\pi^*\omega$ on $M\times \mathbb{R}^2$ is given by $\tilde{\omega}=p_1^*\omega_M+p_2^*\omega_{\text{std}}$, where $p_1:M\times \mathbb{R}^2\rightarrow M$ is the projection to the first factor. Let $\lambda_{\text{std}}$ be a primitive of $\omega_{\text{std}}$ on $\mathbb{R}^2$. Since $p_2\circ u_R:\partial \mathbb{D}\rightarrow \partial \mathbb{D}_R$ is a map of degree one and since $\lambda_{\text{std}}$ restricts to a top degree form on $\partial \mathbb{D}_R\subset \mathbb{R}^2$, we have $$\int_{\partial{\mathbb{D}}}(p_2\circ u_R)^*\lambda_{\text{std}}= \stackrel{=1}{\overbrace{\operatorname{deg}(p_2\circ u_R|_{\partial\mathbb{D}}:\partial \mathbb{D}\rightarrow \partial \mathbb{D}_R)}} \int_{\partial\mathbb{D}_R}\lambda_{\text{std}}=\pi R^2.$$
Moreover, $u_R$ is $\tilde{J}$-holomorphic and thus the area form $dV_{{u_R}^*\tilde{g}}$ of the degenerate conformal metric $u_R^*\tilde{g}$ coincides with $u_R^*\tilde{\omega}$ on $\mathbb{D}$ (with the standard orientation). Using this and $\omega=p_1^*\omega_M+p_2^*\omega_{\text{std}}$, we have
    \begin{align*}
       \text{area}_{u^*\tilde{g}}(\mathbb{D})
        =&
        \int_{\mathbb{D}}{u_R}^*{\widetilde{\omega}}
=\int_{\mathbb{D}}(p_1\circ u_R)^*\omega_M+
        \int_{\mathbb{D}}(p_2\circ u_R)^*\omega_{\text{std}}
        \\
     =&0+\int_{\partial\mathbb{D}}(p_2\circ u_R)^*\lambda_{\text{std}}
        =
        \pi R^2.
    \end{align*}
In passing to the second line we use that
$\int_{\mathbb{D}}(p_1\circ u_R)^*\omega_M=0$ since $\omega_M$ vanishes on the relative homotopy group $\pi_2(M,L)$.
\end{proof}

\begin{proof}[Proof of Theorem~\ref{main; theorem}]
Fix an auxiliary Riemannian metric $g_M$ on $M$. Since $M$ is closed, there is some $K>0$ such that $g_{M}\oplus g_{\mathbb{T}^2}\leq K^2g$,
where \(g_{\mathbb{T}^2}\) denotes a fixed flat metric on \(\mathbb{T}^2\).
This relation lifts to $g_M\oplus g_{\mathbb{R}^2}\leq K^2 \tilde{g}$ where $g_{\mathbb{R}^2}$ denotes the Euclidean metric on $\mathbb{R}^2$. Let $R>0$.\\ By Theorem~\ref{Lem; existence of holomorphic disk},
there is a \(\tilde{J}\)-holomorphic disk $u_R: (\mathbb{D},\partial\mathbb{D})\rightarrow
(M\times \mathbb{R}^2, L\times \partial\mathbb{D}_R)$
such that \(\deg(p_2\circ u_R|_{\partial \mathbb{D}}:\partial \mathbb{D}\rightarrow \partial \mathbb{D}_R)=1\).\\
Now $u_R^*\tilde{g}$ is a (possibly degenerate) Riemannian metric on $\mathbb{D}$ and we can consider the associated (possibly degenerate) distance function $\text{dist}_{u_R^*\tilde{g}}$ on $\mathbb{D}$. \\Let $d_R:=\text{dist}_{u_R^*\tilde{g}}(\cdot, \partial \mathbb{D})$ denote the distance to $\partial \mathbb{D}\subset\mathbb{D}$. It is 1-Lipschitz with respect to $\text{dist}_{u_R^*\tilde{g}}$.
By an argument involving the coarea formula (c.f. Lemma \ref{lem;coarea} with $d=d_R$), there is a compact 1-dimensional submanifold $\Gamma \subset d_R^{-1}([0,R/3K])$ which is homologous to $\partial\mathbb{D}$ in $d_R^{-1}([0,R/3K])$ and such that $$\text{length}_{u_R^* \tilde{g}}(\Gamma)\leq \text{area}_{u_R^*\tilde{g}}(\mathbb{D})/(R/3K)+1\leq 3\pi R K+1,$$ where the last inequality follows from Proposition \ref{remark; area of D with the pull back metric}.\\
Note that $p_2\circ u_R$ maps $d_R^{-1}([0,R/3K])$ to the complement of $ \mathbb{D}_{R/2}$ in $\mathbb{R}^2$, since $g_{\mathbb{R}^2}\leq K^2g$ and since $\partial \mathbb{D}$ is mapped to $\partial \mathbb{D}_R\subset \mathbb{R}^2$.\\
Since $(p_2\circ u_R)_*[\partial\mathbb{D}]=[\partial \mathbb{D}_{R}]$ by the degree assumption in Theorem \ref{Lem; existence of holomorphic disk}, the class $(p_2\circ u_R)_*[\Gamma]=(p_2\circ u_R]_*[\partial \mathbb{D}]\in H_1(\mathbb{R}^2\setminus \mathbb{D}_{R/2})$ is nonzero.
Thus there is a component $\gamma $ of $\Gamma\subset \mathbb{D}$, such that $(p_2\circ u_R)_*[\gamma]$ is nonzero in $H_1(\mathbb{R}^2\setminus \mathbb{D}_{R/2})$. \\
The circle $\gamma$ bounds a disk $\Delta_R=\Delta_R(\gamma)$ in $\mathbb{D}$. It follows that $p_2(u_R(\Delta_R))\supset \mathbb{D}_{R/2}$ since $0\neq (p_2\circ u_R)_*[\gamma]\in H_1(\mathbb{R}^2\setminus \mathbb{D}_{R/2})\cong\ H_1(\mathbb{R}^2\setminus \{w\})$ for all $w\in \mathbb{D}_{R/2}$, and thus $$K^2 \text{area}_{u_R^*\tilde{g}}(\Delta_R)\geq \text{area}_{(p_2\circ u_R)^*g_{\mathbb{R}^2}}(\Delta_R)\geq \text{area}_{g_{\mathbb{R}^2}}(\mathbb{D}_{R/2})=\pi R^2/4.$$
Summarising, we find for each $R>0$ a closed disk $\Delta_R\subset \mathbb{D}$ with smooth boundary $\partial \Delta_R=\gamma$ and
$$  \text{length}_{u_R^*\tilde{g}}(\partial \Delta_R)\leq 3\pi K R+1\text{ and area}_{u_R^*\tilde{g}}(\Delta_R)\geq \pi R^2/4K^2. $$

Since $\partial \Delta_R$ is smooth, the inverse of the Riemann mapping for the interior of $\Delta_R$ extends by \cite[Theorem 5.5]{pommerenke1992boundary} to a smooth map $\varphi_R:\mathbb{D}\rightarrow \Delta_R$ (with complex linear differential). The composition $v_R:=u_R\circ \varphi_R:\mathbb{D}\rightarrow M\times \mathbb{R}^2$ defines a $\tilde{J}$-holomorphic disk with $\text{length}_{v_R^*\tilde{g}}(\partial \mathbb{D})\leq 3\pi KR+1$
and $\text{area}_{v_R^*\tilde{g}}( \mathbb{D})\geq  \pi R^2/4K^2$.\\
\noindent From the $\tilde{J}$-holomorphic disks $v_R:\mathbb{D}\rightarrow M\times \mathbb{R}^2$ we obtain $J=\pi_*\tilde{J}$-holomorphic disks \(
w_R\coloneqq
\pi\circ v_R
:\mathbb{D}\rightarrow M\times \mathbb{T}^2\). Since $\pi$ is locally isometric, these disks $w_R$ satisfy $\text{area}_{w_R^*\tilde{g}}( \mathbb{D})
    =
   \text{area}_{v_R^*\tilde{g}}( \mathbb{D})$ and $
    \text{length}_{w_R^*\tilde{g}}(\partial  \mathbb{D})
    =\text{length}_{v_R^*\tilde{g}}(\partial  \mathbb{D})$.

Therefore the maps $w_R$ satisfy the assumptions of
Theorem~\ref{thm; Duval} below and thus there is a nonconstant $J$-holomorphic plane $w: \mathbb{C}\rightarrow M\times \mathbb{T}^2$.
\end{proof}

\noindent The following is a version of the main Theorem in \cite{Duval}.
\begin{theorem}
    \label{thm; Duval} Let $(X,h) $ be a closed Riemannian manifold and $J$ an almost complex structure on $X$ which preserves $h$.
    Suppose $(w_k:\mathbb{D}\rightarrow X)_{k\in \mathbb{N}}$
    is a sequence of
    \(J\)-holomorphic disks with
    \(\displaystyle
    \lim_{k\to\infty}
    \frac{\text{length}_{w_k^*h}(\partial \mathbb{D})}{\text{area}_{w_k^*h}(\mathbb{D})}
    =0.
    \)\\
     Then there is a nonconstant \(J\)-holomorphic plane $w:\mathbb{C}\rightarrow X$.\hfill \qed
\end{theorem}

\noindent We conclude by proving Lemma \ref{lem;coarea}, which was used in the proof of Theorem \ref{main; theorem}.\\
Let $h$ be a positive semidefinite symmetric 2-tensor, sometimes also called a (possibly) degenerate Riemannian metric, on $\mathbb{D}$.  Associated to $h$, there is a (possibly) degenerate distance function $\text{dist}_h:\mathbb{D}\times \mathbb{D}\rightarrow \mathbb{R}_{\geq 0}$ which is symmetric and satisfies the triangle inequality. \\We denote by $d:\mathbb{D}\rightarrow \mathbb{R}$ the distance to $\partial \mathbb{D}$: $d(x)=\text{dist}_h(x,\partial \mathbb{D})$. By the triangle inequality for $\text{dist}_h$ we have $|d(x)-d(x')|\leq \text{dist}_h(x,x')$; in other words $d$ is 1-Lipschitz with respect to  $\text{dist}_h$.
\begin{lemma}\label{lem;coarea}
Let $b>0$. Then there is for every $\kappa>0$ a compact smooth 1-manifold $\Gamma\subset d^{-1}([0,b])$ such that $[\Gamma]=[\partial \mathbb{D}]\in H_1(d^{-1}([0,b]))$
with $\text{length}_h(\Gamma)\leq \frac{\text{area}_h(\mathbb{D})}{b}+\kappa.$
\end{lemma}

\begin{proof}
In order to construct $\Gamma$, we will perturb $h$ to make a nondegenerate and approximate $d$ by smooth functions and then we apply the coarea formula for smooth functions on Riemannian manifolds.\\
For $\varepsilon>0$ define a smooth Riemannian metric $h_\varepsilon$ by setting $h_\varepsilon:=h+\varepsilon g_{Eu}$, where $g_{Eu}$ denotes the Euclidean metric on $\mathbb{D}$. Since $h_\varepsilon\geq h$ and thus $\text{dist}_{\varepsilon}\geq \text{dist}_h$, the function $d$ is also $1$-Lipschitz with respect to $\text{dist}_{h_\varepsilon}$.
We can use convolution with a bump function to approximate $d$ by functions $d_\delta$, which are smooth in a neighborhood of $d^{-1}([\delta,b-\delta])\subset \mathbb{D}$, $(1+\delta)$-Lipschitz with respect to $\text{dist}_{h_\varepsilon}$ and $\delta$-close to $d$ and such that $d_\delta$ coincides with $d$ near $\partial \mathbb{D}$ and $d^{-1}([b,\infty))$.\\
Applying the coarea formula for smooth maps on Riemannian manifolds, we find $$(1+\delta)\text{area}_{h_\varepsilon}(d_\delta^{-1}([\delta,b-\delta])\geq \int_{d_\delta^{-1}([\delta,b-\delta])} |\nabla d_\delta|_{h_\varepsilon} \,dV_{h_\varepsilon}  = \int_\delta^{b-\delta} \text{length}_{h_\varepsilon}(d_\delta^{-1}(s)) \, ds.$$
It follows that there is  a regular value $s_0\in (\delta,b-\delta)$ of $d_\delta$ such that $$\text{length}_{h_\varepsilon}(d_\delta^{-1}(s_0))\leq(1+\delta) \text{area}_{h_\varepsilon}(\mathbb{D})/(b-2\delta).$$
By construction, this level set is homologous to $\partial\mathbb{D}$ in $d^{-1}([0,b])$. Since $\text{area}_{h_\varepsilon}(\mathbb{D})$ tends to $\text{area}_h(\mathbb{D})$ as $\varepsilon\rightarrow 0$ and since $h\leq h_\varepsilon$, we find for each $\delta>0$ some compact smooth 1-dimensional submanifold $\Gamma_\delta=d_\delta^{-1}(s_\delta)$ with $$\text{length}_h(\Gamma_\delta)\leq \frac{(1+\delta)\text{area}_h(\mathbb{D})}{b-2\delta}+\delta.$$ Thus given $\kappa>0$ we can pick $\Gamma:=\Gamma_\delta$ for $\delta=\delta(\kappa)>0$ sufficiently small.\end{proof}

\section*{Acknowledgements}
We thank David Bechara Senior for bringing the article \cite{cattalani2025ahlfors_classic} to our attention.
\section*{AI statement}
No generative AI tool was used in the research or preparation of this paper.
\bibliographystyle{plain}
\bibliography{references}

\end{document}